\documentclass[12pt,letterpaper]{amsart}

\usepackage{amsfonts,amssymb}
\usepackage{amsmath}

\numberwithin{equation}{section}

\newtheorem{theorem}{Theorem}

\newtheorem{definition}{Definition}

\newtheorem{lemma}{Lemma}

\newcommand{\N}{\mathbb{N}}
\newcommand{\R}{\mathbb{R}}
\newcommand{\B}{\mathbb{B}}
\newcommand{\T}{\mathbb{T}}
\newcommand{\E}{\mathcal{E}}

\newcommand{\Z}{\mathbb{Z}}

\newcommand{\Q}{\mathbb{Q}}

\numberwithin{lemma}{section} \numberwithin{theorem}{section}
\numberwithin{prop}{section} \numberwithin{remark}{section}
\numberwithin{definition}{section} \numberwithin{corollary}{section}

\title[Rigidity of invariant measures]{Rigidity of measures invariant under the action of a multiplicative semigroup of polynomial growth on \(\T\)}
\author{Manfred Einsiedler, Alexander Fish}
\begin{document}

\maketitle
\section{Introduction}


\noindent
In 1967 Furstenberg, in his seminal paper \cite{F}, proved that any closed subset of the torus \( \T = \mathbb{R} / \Z \) which is invariant under the action of multiplication by \( 2 \) and by \( 3 \) (\(\!\!\!\mod 1\)) is either finite or \( \T \).

\noindent For every \( n \in \N \) we denote by \( T_n \)  the following map on \( \T \):
\[
T_n (x) = nx \,  (\!\!\!\!\!\mod 1), \,\, x \in \T.
\]

\noindent The structure of  (\( T_2, T_3 \))-invariant measures on the torus is still  not fully understood.
\vspace{2mm}

\noindent
\textbf{Problem} (Furstenberg): \textit{Is it true that a \textnormal{(\( T_2, T_3 \))}-invariant ergodic
Borel probability  measure on \( \T \) is either Lebesgue or has finite support?}
\vspace{2mm}

\noindent The best known result regarding Furstenberg's measure rigidity problem is due to Rudolph \cite{Rudolph} who proved measure rigidity for the (\( T_2, T_3 \)) action under the additional assumption of positivity of entropy.
\begin{theorem}\textnormal{(Rudolph)}
\label{rudolph_theorem}
Let \( \mu \) be \textnormal{(\( T_2, T_3 \))}-invariant ergodic probability Borel measure on \( \T = \mathbb{R} / \Z \). Then either \( h_{\mu}(T_2) = h_{\mu}(T_3) = 0 \) or \( \mu \) is Lebesgue measure.
\end{theorem}

\noindent
We prove that if a Borel probability measure \( \mu \) on \( \T \) is invariant under the action of a ``large'' multiplicative semigroup (lower logarithmic density is positive) and the action of the whole semigroup is ergodic then \( \mu \)  is either Lebesgue or has finite support, without any entropy assumption.
\begin{theorem}
\label{main_theorem}
Let $\alpha > 0$ and
let $\Sigma \subset \N$ be a multiplicative semigroup with
\begin{equation}
\label{main_equation}
 \liminf_{N \to \infty} \frac{\#\left(\Sigma
\cap [1,N]\right)}{N^{\alpha}} > 0.
\end{equation}
Let $ \mu $ be an ergodic $\Sigma$-invariant Borel probability measure on \( \T = \mathbb{R} / \Z\). Then either \( \mu \) is
supported on a finite number of points or $\mu$ is Lebesgue measure.
\end{theorem}

We note that the assumption to the theorem is equivalent to 
\[
 \liminf_{N \to \infty} \frac{\log\#\left(\Sigma
\cap [1,N]\right)}{\log N} > 0.
\]
Also note that
if \( \Sigma \) has a positive lower density ($\alpha=1$) then  the statement of the theorem follows by an application of the Wiener lemma.


This work is motivated by recent results of Bourgain,
Furman, Lindenstrauss  and Mozes (see \cite{BFLM}) but it is much
simpler. Our case is much easier because the semigroup \( \Sigma \)
is abelian. The largeness of \( \Sigma \) makes it possible to
deduce that zero entropy implies  finite support for \( \mu
\). So, the remaining case is of positive entropy. In this case we
use Johnson's result from \cite{Johnson}, which we state below and is a generalization
of Rudolph's theorem (Thm.~\ref{rudolph_theorem}).

This paper shows, in the concrete setting of endomorphisms of \( \T \),
  that positivity of entropy  is a checkable  condition (see Lemma \ref{lemma2}). The reader is strongly recommended to consult the papers  \cite{BL}, \cite{EKL}, \cite{EL}, \cite{ELMV1}, \cite{ELMV2}, \cite{L} and \cite{MV}, to see   examples of number-theoretic implications  of the checkability of positivity of entropy in more complex situations.

 The authors would like to thank A. Furman, E. Lindenstrauss, F.Nazarov and T. Ward for fruitful discussions.

\section{Main ingredients}\label{sec2}
We remind the reader of the Shannon-McMillan-Breiman theorem.
\newline
\begin{theorem}\label{SMB}
Let \( (X,\B,\mu,T)\) be a measure preserving system such that \( \B \) is a countably generated \(\sigma\)-algebra, and let   \( \xi \) be a countable measurable partition of \( X \) with finite entropy
\textnormal{(\( H_{\mu}(\xi) = \sum_{A \in \xi} -\log(\mu(A)) \mu(A)  < \infty\))}.
Then
\[
\frac{1}{n} I_{\mu}\left( \bigvee_{i=0}^{n-1} T^{-i} (\xi) \right)(x) \to h_{\mu_{x}^{\E}}(T,\xi)
\]
for \( \mu \)-almost every \( x \in X \) and in \(L_{\mu}^1\),  where \( I_{\mu}\left(\xi\right) \) is the information function of the partition \( \xi  \) defined by
\[
I_{\mu}(\xi)(x) = -\log \mu(A), \,\, x \in A \in \xi,
\]
the measures  \( \{ {\mu_{x}^{\E}}\}_{x \in X} \) are the ergodic measures for the transformation \( T \) in the ergodic decomposition of \( \mu \):
\[
 \mu = \int_{X} \mu_{x}^{{\E}} d \mu,
\]
and \( h_{\mu_{x}^{\E}}(T,\xi) \)  denotes  the entropy of \( T\) with respect to \( \xi \) and the measure \( \mu_{x}^{\E} \).
\end{theorem}

\noindent
The other ingredient is Johnson's theorem (\cite{Johnson}) which is a generalization of Rudolph's theorem \ref{rudolph_theorem}.

\noindent We recall the  notion of a nonlacunary multiplicative
semigroup.
\begin{definition}
A multiplicative semigroup \( \Sigma \subset \N \) is called \textbf{nonlacunary} if there does not exist \( a \in \N \) such that \( \Sigma \subset \{ a^n \, | \, n \in \N \cup \{0\}\} \).
\end{definition}
\begin{theorem} \textnormal{(Johnson)}
If \( \Sigma \) is a nonlacunary multiplicative semigroup of integers whose action on \( \T \) as multiplication \( (\!\!\!\mod 1) \) has an ergodic invariant Borel probability measure \( \mu \), then either \( \mu \) is Lebesgue measure or the entropy of each map \( T_n, \, n \in \Sigma \), has \(\mu \)-entropy zero.
\end{theorem}

\section{Proof of Theorem \ref{main_theorem}}

We start with a reformulation of the Shannon-Mcmillan-Breiman theorem (Thm.~\ref{SMB}) in the setting considered here.

\begin{lemma}
\label{lemma1}
In the setting of Theorem \ref{main_theorem},
if \( h_{\mu}(T_p) =0 \) for some \( p \in \Sigma \), then
for any \( \varepsilon > 0 \) and any \( \beta > 0 \), there exist \( \delta_0 > 0 \) and \( X' \subset \T \) of measure \( \mu(X') > 1-\varepsilon \), such that for every positive \( \delta \leq \delta_0 \) and every \( x \in X' \),
\[
\mu(B_{\delta}(x)) > \delta^{\beta},
\]
where \( B_{\delta}(x) \) is the ball of radius \( \delta \) with center at \(x\).
\end{lemma}

One should think of the conclusion as saying that $\mu$ has zero dimension.

\begin{proof} Let \( p \in \Sigma \) such that \( h_{\mu}(T_p) =0 \). Take the partition
\[
\xi = \left\{\left(0,\frac{1}{p}\right],\left(\frac{1}{p},\frac{2}{p}\right],\ldots,\left(\frac{p-1}{p},1\right]\right\}
\]
 of \(\T\). Then
\( \bigvee_{i=0}^{n-1} T_{p}^{-i} (\xi) = \{ (0,\frac{1}{p^n}],\ldots,(\frac{p^n-1}{p^n},1]\}\). Therefore, for every
\( x \in \left(\frac{k}{p^n},\frac{k+1}{p^n}\right]\) and  \( k \in \{0,1,\ldots, p^n -1\} \) we have
\[
I_{\mu}\left( \bigvee_{i=0}^{n-1} T_{p}^{-i} (\xi)\right)(x) = -\log \mu\left(\left(\frac{k}{p^n},\frac{k+1}{p^n}\right]\right).
\]
Let \( \beta > 0 \), without loss of generality we may assume $\beta\in(0,1)$.

 By the a.e.\ convergence in the Shannon-McMillan-Breiman
theorem it follows that for every \( \varepsilon > 0 \) there exists
a measurable set \( X' \subset \T \) with \( \mu(X') > 1-
\varepsilon \) and an \( N \) such that for every \( n \geq N \) we
have
\[
\frac{1}{n}I_{\mu}\left( \bigvee_{i=0}^{n-1} T_{p}^{-i} (\xi)\right)(x) < \beta \log p
\]
for every \( x \in X' \). Take \( \delta_0 = p^{-N} \) and  let \(
\delta \) be any positive number with \( \delta \leq \delta_0 \).
Let \(n \) be the smallest  integer such that \( p^{-n} <
\delta \) (\( n \geq N \)). If \( x \in A =\left(\frac{k}{p^{n}},\frac{k+1}{p^{n}}\right] \in \bigvee_{i=0}^{n-1}
T_{p}^{-i} (\xi)\)  then
 \( A \subset  B_{\delta}(x)  \). Therefore, for every \( x \in X' \), we have
 \[
  I_{\mu}\left( \bigvee_{i=0}^{n-1} T_{p}^{-i} (\xi)\right)(x) =  -\log \mu\left(A\right) \geq -\log \mu(B_{\delta}(x)).
 \]
 So we have, for every \( x \in X' \),
 \[
 \frac{1}{n} \left(-\log \mu(B_{\delta}(x))\right) \leq \frac{1}{n}I_{\mu}\left( \bigvee_{i=0}^{n-1} T_{p}^{-i} (\xi)\right)(x) < \beta \log p.
 \]
 The latter implies that
 \[
 \mu(B_{\delta}(x)) > e^{-n \beta \log p} \geq \delta^{ \beta}p^{-\beta}.
 \]
as $p^{-n+1}\geq \delta$ by choice of $n$. This gives the desired conlusion not quite for the
exponent $\beta$ but for all $\beta'>\beta$ (and $\delta'_0<\delta_0$ depending on $\beta'$).
As $\beta>0$ was arbitrary, the lemma follows.
\end{proof}


In the following we will make use of the assumption of positive logarithmic density of $\Sigma$.

\begin{lemma}
\label{lemma2}
Under the assumptions of Theorem \ref{main_theorem},  if \(h_{\mu}(T_p) = 0 \) for some \( p \in \Sigma \), then the measure \( \mu \)  has  finite support.
\end{lemma}

\begin{proof}
Recall that \( \alpha > 0 \) satisfies the condition
 (\ref{main_equation}). We define $\beta=\frac\alpha{20}$.
Also by \eqref{main_equation} there exists  \( M_0 \) such that for $M\geq M_0$ we have 
 \( \#\left( \Sigma \cap [1,M] \right) > M^{\frac{\alpha}{2}} \).
We also define \( \delta=\delta(M)=M^{-5 } \).

Assume that \( h_{\mu}(T_p) =0 \) for some \( p \in \Sigma \).
From Lemma \ref{lemma1} it follows that for every \( \varepsilon > 0 \) --- $\varepsilon=\frac12$ will do ---  there exists a set \(X' \subset \T \) such that \( \mu(X') > 1 - \varepsilon \)   and  there exists  \( \delta_0 > 0 \) such that for every positive \( \delta < \delta_0 \) and every \( x \in X' \),
\[
\mu(B_{\delta}(x)) > \delta^{\beta}.
\]
We may assume $\delta(M)\leq \delta_0$
for $M\geq M_0$ (by adjusting $M_0$ if necessary). 

At this point we use the invariance\footnote{We will use invariance in the `wrong way' by using that the image of an interval has at least the measure of the original interval.} of \( \mu \) under \( \Sigma \).
Fix some $x\in X'$ and $\delta\leq\delta_0$. For every \( q \in \Sigma \cap [1,M] \) we write \( A_q = q B_{\delta}(x) \). Then \( B_{\delta}(x) \subset T_q^{-1}(A_q) \), and so
by invariance of \( \mu \) under the action of \( T_q \) we get
\[ \mu (A_q) =\mu( q B_{\delta}(x)) =  \mu(T_q^{-1}(A_q)) \geq \mu(B_{\delta}(x)) > \delta^{\beta}. \]
We vary $q\in \Sigma \cap [1,M] $ and note that 
\[
 M^{\frac{\alpha}{2}} \delta^{\beta} = M^{ \frac\alpha{2}-5\frac\alpha{20}} = M^{\frac\alpha{4}} > 1.
\]
To summarize, if we restrict ourselves to $q\in\Sigma\cap[1,M]$ --- and there are at least $M^{\frac\alpha{2}}$ many such choices, every image interval $A_q=q B_\delta(x)$ has $\mu$-measure at least
$\delta^\beta$. In total this would contradict the assumption of having a probability space unless the various sets $A_q$ are not all disjoint.

Hence there exist for any given \( x \in X' \)
\[
 q_1,q_2 \in \Sigma \cap [1,M] , \, q_1 > q_2
 \]
  with
\[
q_1B_{\delta}(x) \cap q_2 B_{\delta}(x) \neq \emptyset.
\]
Thus there exist \( i_1,i_2 \in B_{\delta}(x) \) such that \( q_1 i_1 - q_2 i_2  = k \in \Z \).

We now think of $x\in[0,1)$  and $i_1,i_2\in\R$ as real numbers instead of as cosets belonging to $\mathbb{T}=\R/\Z$. The formulae
\[
k = q_1 ( i_1 - x + x) - q_2 (i_2 - x + x) = q_1 (i_1 - x) - q_2 (i_2 - x) + (q_1 - q_2)x
\]
and
\[
|q_1 (i_1 - x) - q_2 (i_2 - x)| \leq 2 \delta M = 2 M^{-4}
\]
imply that
\[
x = \frac{k}{q_1-q_2} + \kappa, \,\, |\kappa| \leq 2 M^{-4}.
\]
Let  \( \ell = q_1 -q_2 \), so \( \ell < M \).
We summarize: We have shown that $x$ is close to a rational number $\frac{k}\ell$ with related bounds on 
the denominator $\ell$ and on how close $x$ is to $\frac{k}\ell$. Actually the reader should be surprised here on how good a bound on the error $\kappa$ we have achieved --- this will be important. 

For a given $x\in X'$ let us denote the rational number $\frac{k}\ell$ obtained above with the given $M$ as $r(M)=\frac{k}\ell$.  Take some $M_1\geq M_0$. and denote by \( \frac{k_1}{\ell_1}\) the rational \( r(M_1) \) obtained with the above properties. Similarly denote by \( \frac{k_2}{\ell_2}\) the rational \( r(M_1^2) \). So, we have
\[
\left|x - \frac{k_2}{\ell_2} \right| \leq 2 M_1^{-8}
\]
and
\[
\ell_2 \leq M_1^{2}.
\]

\noindent
If we suppose that \( \frac{k_1}{\ell_1} \neq \frac{k_2}{\ell_2} \) then
\[
\left| \frac{k_1}{\ell_1} - \frac{k_2}{\ell_2} \right| \geq \frac{1}{\ell_1 \ell_2} \geq M_1^{-3 } 
\]
On the other hand,
\[
\left|  \frac{k_1}{\ell_1} - \frac{k_2}{\ell_2} \right| \leq \left| \frac{k_1}{\ell_1} - x \right| + \left| \frac{k_2}{\ell_2} - x \right|
\leq 2 (M_1^{-4} + M_1^{-8}).
\]
For \( M_1 \) sufficiently large (there exists positive \( M_0' \geq M_0 \) such that \( M_1 \) should satisfy  \( M_1 \geq M_0' \)), we get a contradiction. Thus \( \frac{k_1}{\ell_1} = \frac{k_2}{\ell_2}\).
Repeating the argument for the given \( x \in X' \) infinitely often (continuing with $M_1^2,M_1^4,\ldots$ replacing $M_1$ in that order) improves the accuracy of the approximation without changing the rational number and so shows that \( x = \frac{k_1}{\ell_1} \).

 The last argument proves that \( X' \) has only finitely many points but has $\mu$-measure at least $\frac12$. So \( \mu \) has atoms. By ergodicity, we get that \( \mu \) has a finite support.
\end{proof}

\noindent
Theorem \ref{main_theorem} follows from A.Johnson theorem cited in Section \ref{sec2} combined with Lemma \ref{lemma2}.

\section{Further Discussion}

 Theorem \ref{main_theorem} might suggest that the reason for the rigidity of invariant measures is some kind of equidistribution phenomenon.

 It was pointed out to us by Fedor Nazarov that for every \( \alpha \not \in \Q \) one can construct a multiplicative semigroup \( \Sigma = \{\sigma_1 < \sigma_2 < \ldots < \sigma_n < \ldots \} \) of positive lower density such that  the sequence \(\{ \sigma_n \alpha \}\) is not equidistributed in \( \T \).
We thank Fedor Nazarov for allowing us to reproduce the construction here.

Our semigroup \( \Sigma \) will satisfy that
\[ d_{*}(\Sigma)  = \liminf_{N \to \infty} \frac{\#\left(\Sigma \cap \{1,\ldots,N\}\right)}{N} \geq \frac{1}{200} \].

 Choose \( \alpha \not \in \Q \). We construct \( \Sigma\) by an iterative process. 

\noindent
Let \( S \subset \N \). For any \( N \in \N \) denote by \( S_N = S \cap \{1,\ldots,N\} \).

\noindent
Let \( N_0 \in \N \) be such that for any \( n \geq N_0 \) we have
\begin{equation}
\label{estimate}
\frac{n}{8} - \frac{n}{1000} < \left|\left\{k  \leq n | 0 < k \alpha (\!\!\!\! \mod 1) < \frac{1}{8}\right\}\right| < \frac{n}{8} + \frac{n}{1000}.
\end{equation}

 We choose $N_1=2N_0$ and let \( B_{N_1} \) to be the subset of size at least \( \frac1{10}N_0 \) of \( \{N_0+1,\ldots,N_1\} \) such that for every \( k \ \in B_{N_1} \) we have \[ 0 < k \alpha (\!\!\!\! \mod 1) < \frac{1}{8}.\] 
If we take the semigroup \( \Sigma_1 \) generated by \( B_{N_1} \) then  there exists  a smallest \( N_2'=N_12^{\ell_1}\) (with $\ell_1\in\N$) such that 
\begin{equation}
\label{ineq_simple}
  |(\Sigma_1)_{2N_2'}| \leq \frac{2N_2'}{100}.
 \end{equation}

 Let \( N_2 = 2N_2' \). We use the lower bound  of (\ref{estimate}) in the interval \( \{1,\ldots,2N_2'\} \), the upper bound of  (\ref{estimate}) in the interval \( \{1,\ldots,N_2'\}\) and the inequality
(\ref{ineq_simple})
to conclude that we can find \( \frac1{10}N_2' \) elements $k$ in
the set \( \{N_2'+1,\ldots,2N_2'\} \setminus \Sigma_1\)  with  \( 0 < k \alpha (\!\!\!\! \mod 1) < \frac{1}{8} \), denote this set by \( A_2 \). We define \( B_{N_2} = (\Sigma_1)_{N_2} \cup A_2\).

 Then we can repeat that process infinitely many times. We will get a sequence of times \( N_1 < N_2 < \ldots < N_m < \ldots \) and the sequence of finite sets \( B_{N_1} \subset B_{N_2} \subset \ldots \) such that for every \( k \in \N \) we have
\begin{equation}
\label{lack_ud}
Re\left(\frac{1}{|B_{N_k}|} \sum_{\sigma \in B_{N_k}} \exp{(2 \pi i \sigma \alpha)}\right) \geq \frac{1}{20} \frac{\sqrt{2}}{2} - \frac{1}{100} > 0.
\end{equation}
Denote by \( \Sigma = \cup_k B_{N_k} \) and notice that by construction $\Sigma_{N_k}=B_{N_k}$ for any $k\geq 1$.  


It is now relatively straight forward to check that \( \Sigma \) is a multiplicative semigroup, 
and that (by the minimality of the $\ell_k$ in the construction) the lower density of $\Sigma$ is at least $\frac{1}{200}$ as claimed.
Clearly (\ref{lack_ud}) implies that \( \Sigma \alpha \) is not uniformly distributed.

By a similar construction, for any choice of countably many  irrational numbers \( \alpha_1, \ldots,\alpha_n, \ldots \) there exists a multiplicative semigroup of a linear growth (of positive lower density) \( \Sigma \) such that for every \( k \in \N \) we have \( \Sigma \alpha_k \) is not uniformly distributed.




\begin{thebibliography}{200}
\bibitem[BFLM]{BFLM} Bourgain, J.; Furman, A.; Lindenstrauss, E.; Mozes, S. Invariant measures and stiffness for non-abelian groups of toral automorphisms.  C. R. Math. Acad. Sci. Paris  344  (2007),  no. 12, 737--742.

\bibitem[BL]{BL} Bourgain, J.; Lindenstrauss, E. Entropy of quantum limits.  Comm. Math. Phys.  233  (2003),  no. 1, 153--171.

\bibitem[EKL]{EKL} Einsiedler, M.; Katok, A.; Lindenstrauss, E. Invariant measures and the set of exceptions to Littlewood's conjecture.  Ann. of Math. (2)  164  (2006),  no. 2, 513--560.

\bibitem[EL]{EL} Einsiedler, M.; Lindenstrauss, E. Diagonalizable flows on locally homogeneous spaces and number theory.  International Congress of Mathematicians. Vol. II,  1731--1759, Eur. Math. Soc., Z\"{u}rich, 2006.

\bibitem[ELMV1]{ELMV1} Einsiedler, M.; Lindenstrauss, E.; Michel, Ph.; Venkatesh, A. The distribution of periodic torus orbits on homogeneous spaces, preprint.

\bibitem[ELMV2]{ELMV2} Einsiedler, M.; Lindenstrauss, E.; Michel, Ph.; Venkatesh, A. Distribution of periodic torus orbits and Duke's theorem for cubic fields, preprint.

\bibitem[F]{F} Furstenberg, H. Disjointness in ergodic theory, minimal sets, and a problem in Diophantine approximation.  Math. Systems Theory  1  (1967), 1--49.

\bibitem[J]{Johnson}  Johnson, Aimee S. A. Measures on the circle invariant under multiplication by a nonlacunary subsemigroup of the integers.  Israel J. Math.  77  (1992),  no. 1-2, 211--240.

\bibitem[L]{L} Lindenstrauss, E. Invariant measures and arithmetic quantum unique ergodicity.  Ann. of Math. (2)  163  (2006),  no. 1, 165--219.

\bibitem[MV]{MV} Michel, P.; Venkatesh, A. Equidistribution, $L$-functions and ergodic theory: on some problems of Yu. Linnik.  International Congress of Mathematicians. Vol. II,  421--457, Eur. Math. Soc., Z\"{u}rich, 2006.

\bibitem[R]{Rudolph} Rudolph, Daniel J. $\times 2$ and $\times 3$ invariant measures and entropy.  Ergodic Theory Dynam. Systems  10  (1990),  no. 2, 395--406.

\end{thebibliography}
\end{document}